\documentclass[11pt, a4paper, twoside]{article}
\usepackage[a4paper]{anysize}
\usepackage{amsmath, amsfonts, amssymb, amsthm, graphicx}
\usepackage[frenchb]{babel}
\usepackage{amsmath,amssymb}
\usepackage[latin1]{inputenc}

\renewcommand{\theta}{\vartheta}

\newcommand{\eps}{\varepsilon}

\newcommand{\A}{\mathbf{A}}

\newcommand{\IP}{\mathbb{P}}

\newcommand{\IR}{\mathbb{R}}
\newcommand{\R}{\mathbb{R}}

\newcommand{\IE}{\mathbb{E}}

\def\y {\mathbf{y}}
\def\x{\mathbf{x}}
\def\w{\mathbf{w}}
\def\b{\mathbf{b}}

\def\v {\mathbf{v}}
\def\f {\mathbf{f}}
\newcommand{\be}{\begin{eqnarray*}}
\newcommand{\ee}{\end{eqnarray*}}
\newcommand{\ben}{\begin{eqnarray}}
\newcommand{\een}{\end{eqnarray}}

\theoremstyle{plain}

\theoremstyle{definition}

\begin{document}

\title{Large deviation for lasso diffusion process}

\author{Azzouz Dermoune $^{1}$,   Khalifa Es-Sebaiy $^{2}$, Youssef Ouknine $ ^{2}$    \vspace*{0.1in} \\
$^{1}$  Laboratoire Paul Painlev\'e, USTL-UMR-CNRS 8524, Lille, France  \\
 Email:azzouz.dermoune@univ-lille1.fr,\\
$^{2}$ Cadi Ayyad University Av. Abdelkrim Khattabi, 40000,
Gu\'eliz-Marrakech,
Morocco \\
Emails: k.essebaiy@uca.ma, ouknine@uca.ma\\
\vspace*{0.1in} } \maketitle {\small \noindent {\bf Abstract:}} The
aim of the present paper is to extend the  large deviation with
discontinuous statistics studied in \cite{BDE} to the     diffusion
$ d\x^\eps=-\{\A^\top(\A\x^\eps-\y)+\mu sgn(\x^\eps)\}dt+\eps d\w$.
The discontinuity of the drift of the diffusion discussed in
\cite{BDE} is equal to the hyperplane $\{\x\in\IR^d:\quad x_1=0\}$,
however, in our case the discontinuity is more complex and is equal
to the set $\{\x\in\IR^d:\quad \prod_{i=1}^dx_i=0\}$.

\section{Introduction}
Let $\y\in\R^n$ be a given vector, $\A$ be a known matrix which maps
the domain $\R^d$ into the domain $\R^n$ and $\mu >0$ is a given
positive real number. The sign of the real number $u$ equals
$sgn(u)=1$ if $u >0$, $sgn(u)=-1$ if $u <0$ and $sgn(0)$ is any
element of $[-1,1]$. The column vector $sgn(\x):=(sgn(x_1), \ldots,
sgn(x_p))^\top$. The following diffusion \ben
d\x^\eps=-\{\A^\top(\A\x^\eps-\y)+\mu sgn(\x^\eps)\}dt+\eps
d\w,\quad \x^\eps(0)=\x(0)\quad \mbox{is given}
\label{lassodiffusion} \een has a discontinuous drift. Using the
fact that $\A^\top(\A\x-\y)+\mu sgn(\x)$ is the subdifferential of
the convex map \be \frac{\|\A\x-\y\|^2}{2}+\mu \|\x\|_1 \ee we can
show that the latter stochastic differential equation (sde) has a
unique strong solution for any $\eps >0$. See (\cite{kree, cepa98,
cepa95, pettersson}). Here $\|\cdot\|$, $\|\cdot\|_1$ denote
respectively the $l2$ and $l1$ norms.

The asymptotic property as $t\to +\infty$ is also possible. The
probability density function \be
\frac{1}{Z}\exp\{-\frac{2}{\eps}(\frac{\|\A\x-\y\|^2}{2}+\mu\|\x\|_1)\}:=p^\eps(d\x)
\ee is the unique invariant probability measure of
$(\frac{\x_t^\eps}{\eps})$, see e.g.  \cite{ABR}. The mode of
$p^\eps$ was introduced in linear regression by \cite{tibshirani96}
and is called lasso. Lasso is the compact and convex set solution of
the system \ben \A_i^\top(\A\x-\y)+\mu sgn(x_i)=0,\quad i=1, \ldots,
d. \label{lasso} \een Here $\A_i^\top$ denotes the $i$-th row of the
matrix $\A^\top$. A large number of theoretical results has been
provided for lasso see e.g.   \cite{CDS, DOR, tibshirani96,
tibshirani13, ounaissi} and the references herein.

If $(P_t^\eps)$ is the semi-groupe defined by (\ref{lassodiffusion})
then  we have the following exponential convergence \ben
\IE_{p^\eps}[|P_t^\eps f-\IE_{p^\eps}(f)|^2]\leq \exp(-t/C)
var_{p^\eps}(f), \label{exponentialconvergence} \een where \be
C=4\IE_{p^\eps}[\|\x-\IE_{p^\eps}(\x)\|^2]. \ee The proof is a
consequence of Poincar\'e inequality (\cite{KLS, bobkov}): \be
var_{p^\eps}(f)\leq
4\IE_{p^\eps}[\|\x-\IE_{p^\eps}(\x)\|^2]\IE_{p^\eps}(\|\nabla\,f\|^2)
\ee valid for all lipschitz map $f$, because $p^\eps$ is
log-concave, and the fact that Poincar\'e inequality is equivalent
to the exponential convergence (\ref{exponentialconvergence}). As a
consequence we can suppose that a.s. $\sup_{t\geq 0}\|\x^\eps(t)\| <
+\infty$.

\section{Limit as $\eps\to 0$}
Let $U:\IR^d\to \IR$ be a convex map such that \be
|\nabla\,U(x)|\leq L(1+|x|),\quad \forall\,x, \ee where $L$ is a
positive constant and $\nabla\,U$ denotes the sub-differential of
$U$.

It's known (see e.g. \cite{kree, cepa98}) that the sde \be
d\x^{\eps}(t)\in-\nabla\,U(\x^\eps(t)) dt+\eps d\w_t, t\in [0,T] \ee
with fixed initial value $\x(0)$,
has a unique strong solution. More precisely there exists
a unique solution
\be \x^\eps(t)=\x(0)-\int_0^t \v^\eps(s)ds+\eps w_t,\quad \forall
t\in [0,T]\ee
where the
measurable map $\v^\eps$ is such that \be \v^\eps(t)\in
\nabla\,U(x_t^\eps)\quad dt\ \mbox{a.e.}\ee

From the linear growth of $\nabla\,U$ we have \be \|\x^\eps(t)\|\leq
K+L\int_0^t |x^\eps(s)|ds, \ee where \be K=\|x(0)\|+LT+\eps \sup_{t\in
[0,T]}\|\w_t\|. \ee Gronwall's lemma tells us that \be \|\x^\eps(t)\|\leq
K\exp(LT), \quad \forall\,t\in [0,T]. \ee Using Ascoli theorem, we
can extract a subsequence such that $\x^\eps\to \x$ uniformely in
$[0,T]$. Now using the inequality \be \|\v^\eps(t)\|\leq
L(1+\|\x^\eps(t)\|)\leq C,\quad \forall\,t, \ee we derive that the
sequence $(\v^\eps)$ is weakly precompact in $L^p([0,T])$ for all
$1< p < +\infty$. Using Mazur's lemma we can construct a measurable
map $\v$ and a subsequence such that $\v^\eps(t)\to \v(t)$, $dt$
a.e.

From the condition $\v^\eps(t)\in \partial\,U(\x^\eps(t))$, the
convergence $(\x^\eps(t),\v^\eps(t))\to (\x(t),\v(t))$ and the fact
that $\partial U$ is monotone maximal we have $v(t)\in
\partial\,U(x(t))$ dt a.e.. Finally the limit $x$ is the unique solution
of the differential inclusion \be dx_t\in -\partial\,U(x_t)dt. \ee
See also \cite{BOQ}.

As an application the solution $\x^\eps$ of (\ref{lassodiffusion})
converges as $\eps\to 0$ to the inclusion equation \ben
d\x^0(t)\in-\{\A^\top(\A\x^0(t)-\y)+\mu
sgn(\x^0)\}dt,\quad\x^0(0)=\x(0). \label{lasso} \een The solution
$\x^0(t)$ converges to lasso as $t\to +\infty$ i.e. $\x^0(t)$
converges to the minimizers of \be
\frac{\|\A\x-\y\|^2}{2}+\mu\|\x\|_1. \ee The set of the latter
minimizers is compact. Hence \ben \sup_{t\geq 0}\|\x^0(t)\| <
+\infty. \label{lassoboundeness} \een Then we have for some $C >0$
and small $\eps$ that \ben \sup_{t\geq 0}\|\x^\eps(t)\| < C
\label{lassoboundeness} \een with a big a probability.

The aim of our work is to extend Boué, Dupuis, Ellis large deviation with discontinuous statistics  \cite{BDE} to the diffusion (\ref{lassodiffusion}).
The discontinuity in \cite{BDE} is equal to the hyperplane
$\{\x\in\IR^d:\quad x_1=0\}$. The discontinuity of the drift of the
diffusion (\ref{lassodiffusion}) is more complex and is equal to the
set $\{\x\in\IR^d:\quad \prod_{i=1}^dx_i=0\}$.

Before arriving to large deviation result we need some preliminary results.
%\end{document}
\section{Preliminary results}
We work in the canonical probability space $(\Omega, \mathcal{F},\IP)$
where $\Omega=C([0,1],\IR^d)$ endowed with its Borel $\sigma$-field $\mathcal{F}$,
and its Wiener measure $\IP$.
The canonical process $W_t: \w\in \Omega\to \w(t)$, $t\in [0,1]$ is the Wiener process
under $\IP$. The filtration $\mathcal{F}_t:=\sigma(\{W_s:\quad s\leq t\}, \mathcal{N})$,
$t\in [0,1]$,
where $\mathcal{N}$ is the collection of the $\IP$-null sets.
The diffusion $\x^\eps$ (\ref{lassodiffusion}) is considered in the canonical
probability space $(\Omega, \mathcal{F},\IP)$.
Its discontinuous drift is
\ben
\b(\x):=-\{\A^\top(\A\x-\y)+\mu sgn(\x)\}.
\label{lassodrift}
\een
 We denote by
$\IE_{\x(0)}$ the mathematical expectation under
the probability distribution of the solution $\x^\eps$
known that $\x^\eps(0)=\x(0)$.

I) We define for each $i=1, \ldots, d$, the
Borel measures \be
\gamma_i^{\eps,1}(dt)={\bf 1}_{[x_i^\eps(t)\leq 0]}dt,\\
\hat{\gamma}_i^{\eps,1}(t)={\bf 1}_{[x_i^\eps(t)\leq 0]},\\
\gamma_i^{\eps,2}(dt)={\bf 1}_{[x_i^\eps(t)> 0]}dt,\\
\hat{\gamma}_i^{\eps,2}(t)={\bf 1}_{[x_i^\eps(t)> 0]}. \ee By
extracting a subsequence we have $\x^\eps\to \x^0$ where $\x^0$ is
the solution of the inclusion differential equation (\ref{lasso}),
and \be (\gamma_i^{\eps,\eta}(dt),i=1, \ldots, d, \eta=1, 2)\to
(\gamma_i^{\eta}(dt),i=1, \ldots, d, \eta=1, 2) \ee where the Borel
measures $(\gamma_i^{\eta}(dt),i=1, \ldots, d, \eta=1, 2)$ satisfy
\ben
&&\gamma_i^{\eta}(dt)=\hat{\gamma}_i^\eta(t)dt,\forall\,i=1, \ldots, d, \eta=1, 2,\nonumber\\
&&\hat{\gamma}_i^1(t)+\hat{\gamma}_i^2(t)=1, \forall\,i=1, \ldots, d,\nonumber\\
&&\hat{\gamma}_i^1(t)=1, \quad \mbox{if} \,\, x_i^0(t) < 0,\nonumber\\
&&\hat{\gamma}_i^2(t)=1, \quad \mbox{if} \,\, x_i^0(t) > 0,\nonumber\\
&&\hat{\gamma}_i^2(t)-\hat{\gamma}_i^1(t):=sgn(x_i^0(t))\quad \mbox{if} \,\, x_i^0(t)=0,\nonumber\\
&&-\{\A_i^\top(\A\x^0(t)-\y)-\mu\}\geq 0,\,\,\mbox{and}\,\,
-\{\A_i^\top(\A\x^0(t)-\y)+\mu\}\leq 0\quad \mbox{if} \,\,
x_i^0(t)=0\label{Dupuis}. \een

The property (\ref{Dupuis}) tells us
that $x_i^0(t)$ stays at zero when the strength \be
|\A_i^\top(\A\x^0(t)-\y)|\leq \mu. \ee This phenomenon is known by
physicist \cite{Baule}, and we can show it mathematically  using a similar proof
as in \cite{BDE}.

II) Now we fix $\f$ deterministic such that $\int_0^1\|\f(t)\|dt <
+\infty$. We consider the sde \ben d\x^\eps(t)=\{\f(t)-\mu
sgn(\x^\eps(t))\}dt+\eps d\w(t),\quad \x(0)\quad \mbox{is given},
\label{xf} \een and its limit $\x^0$ as $\eps\to 0$ is the solution
of the differential inclusion \ben d\x^0(t)\in \{\f(t)-\mu
sgn(\x^0(t))\}dt,\quad \x(0)\quad \mbox{is given}. \label{xzero}
\een We have $dt$ a.e. \be
\frac{d\x^{0}(t)}{dt}=\f(t)-\mu\{\hat{\gamma}^2(t)-\hat{\gamma}^1(t)\}.
\ee

1) If $x_i^{0}(t)<0$, then $\hat{\gamma}_i^2(t)=0$, $\hat{\gamma}_i^1(t)=1$ and
\be
\frac{dx_i^{0}(t)}{dt}=f_i(t)+\mu.
\ee

2) If $x_i^{0}(t)>0$, then $\hat{\gamma}_i^1(t)=0$, $\hat{\gamma}_i^2(t)=1$ and
\be
\frac{dx_i^{0}(t)}{dt}=f_i(t)-\mu.
\ee

3) We have $dt$ a.e on the set $\{t:\quad x_i^{0}(t)=0\}$ that
\be
-\mu\leq f_i(t)\leq \mu
\ee
and
\be
\frac{dx_i^{0}(t)}{dt}&=&f_i(t)-\mu\{\hat{\gamma}_i^{2}(t)-\hat{\gamma}_i^{1}(t)\}\\
&=&\hat{\gamma}_i^2(t)\{f_i(t)-\mu\}+\hat{\gamma}_i^1(t)\{f_i(t)+\mu\}\\
&=&0.
\ee
It follows that
\be
f_i(t)&=&\mu\{\hat{\gamma}_i^2(t)-\hat{\gamma}_i^1(t)\},\\
1&=&\hat{\gamma}_i^2(t)+\hat{\gamma}_i^1(t).
\ee
Hence
\be
\hat{\gamma}_i^2(t)=\frac{f_i(t)+\mu}{2\mu},\\
\hat{\gamma}_i^1(t)=\frac{\mu-f_i(t)}{2\mu}. \ee Finally, if
$x_i^{0}(t)=0$, then $dt$ a.e. $\frac{dx_i^{0}(t)}{dt}=0$ and \be
\hat{\gamma}_i^2(t)=\frac{f_i(t)+\mu}{2\mu},\\
\hat{\gamma}_i^1(t)=\frac{\mu-f_i(t)}{2\mu}.
\ee
Observe that $\beta^{1}(t):=f_i(t)+\mu\geq 0$, $\beta^{2}(t):=f_i(t)-\mu\leq 0$, and
$\frac{dx_i^{0}(t)}{dt}=\hat{\gamma}_i^2(t)\beta^{2}(t)+\hat{\gamma}_i^1(t)\beta^{1}(t)$.
Observe also that $x_i^0(t)\neq 0$ if and only if $|f_i(t)| >\mu$.

By choosing $f_i$ piecewise constant we obtain the trajectorie $\x$
having the following properties: There exist $0=\tau_1 <\ldots < \tau_r=1$ and
the constants $(\beta_i(k), i=1, \ldots, d, k=1, \ldots, r)$ such that
\be
&&1) \frac{dx_i(t)}{dt}=\beta_i(k),\quad \forall\,t\in [\tau_k,\tau_{k+1}),\\
&&2) x_i(t)\neq 0,\quad \forall\,t\in [\tau_k,\tau_{k+1}),\quad\mbox{or}\quad
x_i(t)= 0,\quad \forall\,t\in [\tau_k,\tau_{k+1}).
\ee
We denote by $\mathcal{N}_0$ the set of the maps $\x: [0,1]\to \IR^d$
which satisfy the latter properties. It's a dense subset of $C([0,1])$.

III) Now we introduce the set
\be
\mathcal{A}=\{\v:\Omega\times [0,1]\to \IR:\quad \mbox{is progressively measurable and}\quad
\IE_{\x(0)}[\int_0^1\|\v(t)\|^2dt] < +\infty\},
\ee
and for $\v\in\mathcal{A}$ we denote  by $\x^{\eps,\v}$ the solution
\be
d\x^{\eps,\v}(t)=\{\b(\x^{\eps,\v}(t))+\v(t)\}dt+\eps d\w(t),\quad \x^{\eps,\v}(0)=\x(0),
\ee
where
\be
\b(\x)=-\{\A^\top(\A\x^\eps-\y)+\mu sgn(\x)\}.
\ee
Let $(\v^\eps, \eps\in (0,1])\subset \mathcal{A}$ be a family of
progressively measurable processes such that
\ben
\sup_{\eps\in (0,1]}\IE_{\x(0)}[\int_0^1\|\v^\eps(t)\|^2dt] < +\infty.
\label{UI}
\een

We define for each $i=1, \ldots, d$, the Borel measures \be
&&\nu^\eps(d\v,t)dt=\delta_{\v^\eps(t)}(d\v)dt,\\
&&\nu_i^{\eps,1}(d\v,t)={\bf 1}_{[x_i^\eps(t)\leq 0]}\delta_{\v^\eps(t)}(d\v),\\
&&\nu_i^{\eps,2}(d\v,t)={\bf 1}_{[x_i^\eps(t)> 0]}\delta_{\v^\eps(t)}(d\v),\\
&&\gamma_i^{\eps,1}(dt)={\bf 1}_{[x_i^\eps(t)\leq 0]}dt,\\
&&\hat{\gamma}_i^{\eps,1}(t)={\bf 1}_{[x_i^\eps(t)\leq 0]},\\
&&\gamma_i^{\eps,2}(dt)={\bf 1}_{[x_i^\eps(t)> 0]}dt,\\
&&\hat{\gamma}_i^{\eps,2}(t)={\bf 1}_{[x_i^\eps(t)> 0]}. \ee By
extracting a subsequence we have $\x^{\eps,\v^\eps}\to \x^{0}$, and
$\nu^\eps\to \nu$. Thanks to (\ref{UI}), we have (see \cite{BDE})
\be \IE_{\x(0)}[\int_{\IR^d}\|\v\|\nu(d\v,t)] < +\infty. \ee The
limit $\x^{0}$ is the solution of the differential inclusion \be
d\x^{0}(t)\in \{\f(t)-\mu sgn(\x^{0}(t))\}dt, \ee where \be
\f(t)=-\A^\top(\A\x^{0}(t)-\y)+\int_{\IR^d}\v\nu(d\v,t). \ee Hence
$\x^{0}$ is exactly the solution studied in (\ref{xzero}).

%\end{document}
\section{Boué Dupuis variational representation}
The variational representation of
 \cite{BD}  tells us that for any bounded measurable map $h:
(\Omega, \mathcal{F},\IP)\to \IR$ \ben
H^\eps(\x(0)):&=&-\eps^2\ln\left(\IE_{\x(0)}[\exp\{-\frac{h(\x^\eps)}{\eps^2}\}]\right)\\
&=&\inf\{\IE_{\x(0)}[\frac{1}{2}\int_0^1\|\v(t)\|^2dt+h(\x^{\eps,\v})]:\quad
\v\in \mathcal{A}\},
\label{VR}
\een
where $\x^{\eps,\v}$ denotes the solution
\be
d\x^{\eps,\v}(t)=\{\b(\x^{\eps,\v}(t))+\v(t)\}dt+\eps d\w(t),\quad \x^{\eps,\v}(0)=\x(0).
\ee

The control $\v^\eps\in \mathcal{A}$ such that \be H^\eps(\x(0))\geq
\IE_{\x(0)}[\frac{1}{2}\int_0^1\|\v^\eps(t)\|^2dt+h(\x^{\eps,\v^\eps})]-\eps^2
\ee and the diffusion \be
d\x^{\eps,\v^\eps}(t)=\{\b(\x^{\eps,\v^\eps}(t))+\v^\eps(t)\}dt+\eps
d\w(t),\quad \x^{\eps,\v^\eps}(0)=\x(0) \ee play the central role in
the large deviation result \cite{BDE}, and then also in our case. We
set \be \bar{\x}^\eps=\x^{\eps,\v^\eps},\quad \bar{\x}=\lim_{\eps\to
0}\bar{\x}^\eps. \ee Observe  that Condition 3.2. in \cite{BDE} \be
\sup_{\eps\in (0,1]}\IE_{\x(0)}[\int_0^1\|\v^\eps(t)\|^2dt] <
+\infty \ee holds also in our case.

\section{Large deviation: upper bound}
We start from the variational representation (\ref{VR}):
\be
H^\eps(\x(0)):=-\eps^2\ln\{\IE_{\x_0}[\exp(-\frac{h(\x^\eps)}{\eps^2}]\}=\inf_{\v\in\mathcal{A}}\IE_{\x_0}[
\frac{1}{2}\int_0^1\|\v(t)\|^2dt+h(\x^{\eps,v})]
\ee
valid for all bounded measurable map $h$.

From the definition of $\v^\eps$ we have
\be
H^\eps(\x(0)) \geq \IE_{\x_0}[
\frac{1}{2}\int_0^1\|\v^\eps (t)\|^2dt+h(\bar{\x}^{\eps})]-\eps^2.
\ee
It follows that
\be
&&\lim\inf_{\eps\to 0} H^\eps(\x(0)) \geq \lim\inf_{\eps\to 0}\IE_{\x_0}[\frac{1}{2}
\int_{[0,1]\times \IR^d}\|\v\|^2\nu^\eps (d\v,t)dt+h(\bar{\x}^\eps)]\\
&&=\lim\inf_{\eps\to 0}\IE_{\x_0}[\frac{1}{2}\int_{[0,1]\times \IR^d}\|\v\|^2\nu^\eps (d\v,t)dt+h(\bar{\x})].
\ee
From Fatou lemma we have
\be
\lim\inf_{\eps\to 0}\IE_{\x_0}[
\int_{[0,1]\times \IR^d}\|\v\|^2\nu^\eps (d\v,t)dt=
\IE_{\x_0}[\int_0^1\lim\inf_{\eps\to 0}\int_{\IR^d}\|\v\|^2\nu^\eps(d\v,t)dt].
\ee
Using the inequality
\be
\lim\inf_{\eps\to 0}\int f(\v)\mu^n(d\v)\geq \int f(\v)\mu(d\v)
\ee
valid for all $f\geq 0$ measurable and all $\mu^n\to \mu$
weakly, we obtain
\be
\lim\inf_{\eps\to 0}\int_{\IR^d}\|\v\|^2\nu^\eps(d\v,t)\geq \int_{\IR^d}\|\v\|^2\nu(d\v,t).
\ee
Finally we have
\be
\lim\inf_{\eps\to 0} H^\eps(\x(0))\geq
\IE_{\x(0)}[\frac{1}{2}\sum_{i=1}^d\int_{[0,1]\times \IR^d}|v_i|^2\nu(d\v,t)dt+h(\bar{\x})].
\ee
If $\bar{x}_i(t)< 0$, then
\be
\frac{d\bar{x}_i}{dt}(t)=b_i^1(\bar{\x}(t))+\int_{\IR^d}v_i\nu_i^1(d\v,t).
\ee
If $\bar{x}_i(t)> 0$, then
\be
\frac{d\bar{x}_i}{dt}(t)=b_i^2(\bar{\x}(t))+\int_{\IR^d}v_i\nu_i^2(d\v,t).
\ee
We also recall that in these cases
\be
\nu_i^\eta(d\v,t)
\ee
is a probability measure for $\eta=1, 2$. It follows from Jensen inequality that
\be
\int_{\IR^d}|v_i|^2\nu_i^\eta(d\v,t)&\geq &(\int_{\IR^d}v_i\nu_i^\eta(d\v,t))^2\\
&\geq & |\frac{d\bar{x}_i(t)}{dt}-b_i^\eta(\bar{\x}(t))|^2:=L_i^\eta(\bar{\x}(t),
\frac{d\bar{x}_i(t)}{dt}).
\ee
If $\bar{x}_i(t)=0$, then $0<\hat{\gamma}_i^1(t)< 1$, and \be
\nu(d\v,t)&=&\hat{\gamma}_i^1(t)\frac{\nu_i^1(d\v,t)}{\hat{\gamma}_i^1(t)}+
\hat{\gamma}_i^2(t)\frac{\nu_i^2(d\v,t)}{\hat{\gamma}_i^2(t)}. \ee
The measure $\frac{\nu_i^\eta(d\v,t)}{\hat{\gamma}_i^\eta(t)}$ is a
probability for each $\eta=1, 2$. Again from Jensen inequality we have \be
\int_{\IR^d}|v_i|^2\frac{\nu_i^\eta(d\v,t)}{\hat{\gamma}_i^\eta(t)}\geq
\left(\int_{\IR^d}v_i\frac{\nu_i^\eta(d\v,t)}{\hat{\gamma}_i^\eta(t)}\right)^2.
\ee We recall that if $\bar{x}_i(t)=0$, then \be
\beta_i^1(t)&=&b_i^1(\bar{\x}(t))+\int_{\IR^d}v_i\frac{\nu_i^1(d\v,t)}{\hat{\gamma}_i^1(t)}\geq 0,\\
\beta_i^2(t)&=&b_i^2(\bar{\x}(t))+\int_{\IR^d}v_i\frac{\nu_i^2(d\v,t)}{\hat{\gamma}_i^2(t)}\leq 0,
\ee
and if we denote
\be
\beta_i=\frac{d\bar{x}_i}{dt}(t)
\ee
then
\be
\hat{\gamma}_i^1(t)\beta_i^1(t)+\hat{\gamma}_i^2(t)\beta_i^2(t)=\beta_i.
\ee
It follows that
\be
&&\int_{\IR^d}|v_i|^2\nu(d\v,t)\geq \hat{\gamma}_i^1(t)|\beta_i^1(t)-
b_i^1(\bar{\x}(t))|^2+\hat{\gamma}_i^2(t)|\beta_i^2(t)-
b_i^2(\bar{\x}(t))|^2 \geq \\
&&\inf\{p^1|\beta_i^1-b_i^1(\bar{\x}(t))|^2+p^2|\beta_i^2-
b_i^2(\bar{\x}(t))|^2\}:=L_i^0(\bar{\x}(t),\frac{d\bar{x}_i}{dt}(t)).
\ee
The infimum is taken under the constraint
\be
p^1, p^2> 0,\quad p^1+p^2=1,\\
p^1\beta_i^1+p^2\beta_i^2=\frac{d\bar{x}_i}{dt}(t):=\beta_i.
\ee
%\end{document}

We define
\be
\tilde{L}_i(\bar{\x}(t),\frac{d\bar{x}_i}{dt}(t))=|
\frac{d\bar{x}_i}{dt}(t)-b_i^1(\bar{\x}(t))|^2:=L_i^{(1)}(\bar{\x}(t),
\frac{d\bar{x}_i}{dt}(t)),\quad \mbox{if}\quad \bar{x}_i(t) < 0,\\
\tilde{L}_i(\bar{\x}(t),\frac{d\bar{x}_i}{dt}(t))=|\frac{d\bar{x}_i}{dt}(t)-b_i^2(\x(t))|^2:=
L_i^{(2)}(\bar{\x}(t),\frac{d\bar{x}_i}{dt}(t)),\quad \mbox{if}\quad \bar{x}_i(t) > 0,\\
\tilde{L}_i(\bar{\x}(t),\frac{d\bar{x}_i}{dt}(t)):=L_i^0(\bar{\x}(t),\frac{d\bar{x}_i}{dt}(t)),
\quad \mbox{if}\quad \bar{x}_i(t)=0.
\ee
%\end{document}
It follows for each $i$ that \be
\left(\int_{\IR^d}v_i\nu(d\v,t)\right)^2\geq
\tilde{L}_i(\bar{\x}(t),
\frac{d\bar{x}_i}{dt}(t)), \ee and \be
\lim\inf_{\eps\to 0} H^\eps(\x(0))&\geq &\IE_{\x(0)}[
\frac{1}{2}\sum_{i=1}^d\int_{[0,1]}\tilde{L}_i(\bar{\x}(t),
\frac{d\bar{x}_i}{dt}(t))dt+h(\bar{\x})]\\
&\geq &\inf\{\frac{1}{2}I(\varphi)+h(\varphi):\varphi\in C_{\x(0)}([0,1])\},
\ee
where the rate function
\ben
I_{\x(0)}(\varphi):=\sum_{i=1}^d\int_0^1\tilde{L}_i(\varphi(t),\frac{d\varphi_i}{dt}(t))dt.
\label{ratefunction}
\een
%\end{document}
The infimum is equal to $+\infty$ if the latter set is empty.

Finally we have for any sequence $\eps$ such that
$\bar{\x}^\eps\to \bar{\x}$,
and $(\nu^\eps, \nu_i^{\eps,\eta}, \gamma_i^{\eps,\eta}, i=1, \ldots, d, \eta)\to
(\nu, \nu_i^{\eta}, \gamma_i^{\eta},i=1, \ldots, d, \eta)$
that
\be
\lim\inf_{\eps\to 0} H^\eps(\x(0))\geq \inf_{\varphi\in C([0,1])}\{\frac{1}{2}I(\varphi)+h(\varphi)\}.
\ee

Using the same argument as in Boué-Dupuis-Ellis \cite{BDE} we can show that
\be
\lim\inf_{\eps\to 0} H^\eps(\x(0))\geq \inf_{\varphi\in C([0,1])}\{I_{\x_0}(\varphi)+h(\varphi)\}
\ee
for all $\eps\in (0,1]$.
%\end{document}
\section{Large deviation: lower bound}

\subsection{Properties of $L_i^{\eta}, \eta=0,1,2$}
We define $L_i^0:\IR^d\times \IR\to [0,+\infty)$ by
\ben
L_i^0(\x,\beta_i)=\inf\{p^1|\beta_i^1-b_i^1(\x)|^2+p^2|\beta_i^2-
b_i^2(\x)|^2\}
\label{Lzero}
\een
The infimum is taken under the constraint
\be
p^1, p^2> 0,\quad p^1+p^2=1,\quad \beta_i^1\geq 0,\quad \beta_i^2\leq 0\\
p^1\beta_i^1+p^2\beta_i^2=\beta_i.
\ee

{\it {\bf Proposition.} 1) If $(\x,\beta_i)\in\IR^d\times \IR$ are such that
\be
b_i^2(\x)<\beta_i < b_i^1(\x)
\ee
then
\be
L_i^0(\x,\beta_i)=0.
\ee

2) If $\beta_i\leq b_i^2(\x)$ then
\be
L_i^0(\x,\beta_i)=|\beta_i-b_i^2(\x)|^2.
\ee

3) If $\beta_i\geq b_i^2(\x)$ then
\be
L_i^0(\x,\beta_i)=|\beta_i-b_i^1(\x)|^2.
\ee
}
{\bf Proof.} Observe that
\be
p|\beta_i^1-b_i^1(\x)|^2+(1-p)|\beta_i^2-b_i^2(\x)|^2\geq |pb_i^1(\x)+(1-p)b_i^2(\x)-\beta|^2
\ee
for any $p\in (0,1)$ and any couple $\beta_i^1,\beta_i^2$ such that
$p\beta_i^1+(1-p)\beta_i^2=\beta$.
If the infimum (\ref{Lzero}) is such that
\be
L_i^0(\x,\beta)=p|\beta_i^1-b_i^1(\x)|^2+(1-p)|\beta_i^2-b_i^2(\x)|^2
\ee
for some $p\in (0,1)$, $\beta_i^1\geq 0$, $\beta_i^2\leq 0$, then
\be
|\beta_i^1-b_i^1(\x)|^2=|\beta_i^2-b_i^2(\x)|^2=|pb_i^1(\x)+(1-p)b_i^2(\x)-\beta|^2.
\ee
Hence
\be
L_i^0(\x,\beta)=\inf\{|pb_i^1(\x)+(1-p)b_i^2(\x)-\beta|^2\}
\ee
where the infimum is also under the same constraint as in (\ref{Lzero}).
This finishes the proof.
%\end{document}

{\it
{\bf Corollary.} 1) For each $i=1, \ldots, d$, $\eta=0, 1, 2$,
the maps $(\x,\beta)\in\IR^d\times\IR\to L_i^\eta(\x,\beta)$
are continuous.

2) If $\beta_i\leq 0$, then $L_i^0(\x,\beta_i)\leq L_i^2(\x,\beta_i)$.

3) If $\beta_i\geq 0$, then $L_i^0(\x,\beta_i)\leq L_i^1(\x,\beta_i)$.

4) For each $i=1, \ldots, d$ and for $\x$ fixed, the map
$\beta_i\to L_i^0(\x,\beta_i)$ is convex.
}
%\end{document}
%\end{document}
%\end{document}

\bigskip

Now back to the large deviation lower bound.
We are going to show that
\be
\lim\sup_{\eps\to 0} H^\eps(\x(0))\leq I_{\x_0}(\varphi)+h(\varphi)
\ee
for all $\varphi\in \mathcal{N}_0$. The map $\varphi$ is defined
by $(t_k,\beta(k)): k=1, \ldots, r$ such that
\be
\frac{d\varphi_i}{dt}(t)=\beta_i(k),\,t\in [t_k,t_{k+1}[,\\
\varphi_i(t)\neq 0,\,t\in [t_k,t_{k+1}[, \,\mbox{or}\,\,
\varphi_i(t)=0,\,t\in [t_k,t_{k+1}[.
\ee
If $\varphi_i(t)=0$ on $[t_k,t_{k+1}[$, then $\beta_i(k):=\frac{d\varphi_i}{dt}(t)=0$
on $[t_k,t_{k+1}[$.

%\end{document}
We consider the control
\be
\v^1(\x,t)=\beta^1(k)-\b^1(\x),\quad \v^2(\x,t)=\beta^2(k)-\b^2(\x)
\ee
for $t\in [t_k,t_{k+1})$. Here $\beta_i^1(k)=-\mu$, $\beta_i^2(k)=\mu$
if $\beta_i(k)=0$. Observe that $0=\frac{\beta_i^1(k)+\beta_i^2(k)}{2}$.
If $\beta_i(k)\neq 0$, then
$\beta_i^1(k)=\beta_i^2(k)=\beta_i(k)$.

Now define for $i=1, \ldots, d$,
\be
v_i(\x,t)=v_i^1(\x,t){\bf 1}_{[x_i\leq 0]}+v_i^2(\x,t){\bf 1}_{[x_i>0]},
\ee
and $\v(\x,t)$ denotes the vector column with the components $v_i(\x,t)$.
%\end{document}
The controlled process
\ben
&&dx_i^{\eps,\varphi}(t)=b_i(\x^{\eps}(t))dt+v_i(\x^{\eps}(t),t)dt+\eps dw_i(t)\nonumber\\
&&=\{\beta_i^2(k){\bf 1}_{[x_i^{\eps,\varphi}(t)> 0]}+\beta_i^1(k)
{\bf 1}_{[x_i^{\eps,\varphi}(t)\leq 0]}\}dt+\eps dw_i(t)\label{f-mu}
\een
Let us define
\be
f_i(k)=\beta_i(k)+\mu,\quad t\in [t_k,t_{k+1}),\quad \varphi_i(t_k)>0,\\
f_i(k)=\beta_i(k)-\mu,\quad t\in [t_k,t_{k+1}),\quad \varphi_i(t_k)<0,\\
f_i(k)=0,\quad t\in [t_k,t_{k+1}),\quad \varphi_i(t_k)=0. \ee Then
we can rewrite (\ref{f-mu}) as \be
dx_i^{\eps,\varphi}(t)=\{f_i(t)-\mu
sgn(x_i^{\eps,\varphi}(t))\}dt+\eps dw_i(t) \ee where \be
f_i(t)=f_i(k),\quad t\in [t_k,t_{k+1}). \ee It follows from
(\ref{xf}) that $\x^{\eps,\varphi}$ converges to $\varphi$ as
$\eps\to 0$.

From the variational representation
\be
H^\eps(\x(0))=\inf_{\v\in\mathcal{A}}\IE_{\x(0)}[\frac{1}{2}\int_0^1\|\v(t)\|^2dt+h(\x^{\eps,\v})]
\ee
we have
\be
&&\limsup_{\eps\to 0} H^\eps(\x(0))\leq \limsup_{\eps\to 0}
\IE_{\x(0)}[\frac{1}{2}\int_0^1\|\frac{d\varphi}{dt}(t)-\b(\x^{\eps,\varphi}(t))\|^2dt+h(\x^{\eps,\varphi})]\\
&&=[\frac{1}{2}\int_0^1\|\frac{d\varphi}{dt}(t)-\b(\varphi(t))\|^2dt+h(\varphi)],
\ee
where
\be
\int_0^1\|\frac{d\varphi}{dt}(t)-\b(\varphi(t))\|^2dt=\sum_{i=1}^d
\int_0^1|\frac{d\varphi_i}{dt}(t)-b_i(\varphi(t))|^2dt
\ee
and
\be
b_i(\varphi(t))=b_i^{1}(\varphi(t)),\quad \mbox{if}\quad \varphi_i(t) < 0,\\
b_i(\varphi(t))=b_i^{2}(\varphi(t)),\quad \mbox{if}\quad \varphi_i(t) > 0,
\ee
in these cases
\be
|\frac{d\varphi_i}{dt}(t)-b_i(\varphi(t))|^2=\tilde{L}_i(\varphi(t),\frac{d\varphi_i}{dt}(t)).
\ee
Moreover, on each interval $[t_k,t_{k+1})$ such that
$\beta_i(k)=0$ we can show that
\be
|\frac{d\varphi_i}{dt}(t)-b_i(\varphi(t))|^2:=|b_i^1(\beta(k))|^2\quad \mbox{if}\quad b_i^1(\beta(k))\leq 0,\\
|\frac{d\varphi_i}{dt}(t)-b_i(\varphi(t))|^2:=|b_i^2(\beta(k))|^2\quad \mbox{if}\quad b_i^2(\beta(k))\geq 0,\\
|\frac{d\varphi_i}{dt}(t)-b_i(\varphi(t))|^2:=0,\quad b_i^2(\beta(k))< 0 < b_i^1(\beta(k)).
\ee
More precisely, if $\beta_i(k)=0$ then
\be
|\frac{d\varphi_i}{dt}(t)-b_i(\varphi(t))|^2=L_i^0(\varphi(t),\frac{d\varphi_i}{dt}(t)).
\ee
Finally we have
\be
\int_0^1\|\frac{d\varphi}{dt}(t)-\b(\varphi(t))\|^2dt=\int_0^1\tilde{L}(\varphi(t),\frac{d\varphi}{dt}(t))dt
\ee
and then for all $\varphi\in \mathcal{N}_0$
\be
\limsup_{\eps\to 0} H^\eps(\x(0))\leq I_{\x(0)}(\varphi)+h(\varphi).
\ee

To finish the proof of the large deviation's lower bound we need the following
lemmas.
The proof is the same as in Dupuis and Ellis book \cite{DE}. For the sake of completeness
we will recall the proof.

{\bf Lemma 1.}
Let $\psi\in C_{\x(0)}([0,1])$ such that $\int_0^1\tilde{L}(\psi(t),\psi'(t))dt < +\infty$.
For each $\delta >0$, there exist $\theta>0$ and $\xi\in C_{\x(0)}([0,1])$ with the following properties
\be
\sup_{t\in [0,1]}|\xi(t)-\psi(t)| \leq \delta,\quad
\int_0^1\tilde{L}(\xi(t),\frac{d\xi}{dt}(t))dt\leq \int_0^1\tilde{L}(\psi(t),\frac{d\psi}{dt}(t))dt+\delta
\ee
and
\be
\sup_{t\in [0,1]}\|\frac{d\xi}{dt}(t)\| \leq \theta.
\ee

Now we prove the following result.

{\bf proof.} Let $c, \lambda\in (0,1)$. We define
\be
D_\lambda&=&\{t\in [0,1]:\quad \|\frac{d\psi}{dt}(t)\|\geq \frac{1}{\lambda}\},\\
E_\lambda &=&\{t\in [0,1]:\quad \|\frac{d\psi}{dt}(t)\|<\frac{1}{\lambda}\}.
\ee
We construct the time-rescaling map $S_\lambda: [0,1]\to [0,+\infty)$ as follows
\be
S_\lambda(0)=0,\quad \frac{dS_\lambda}{dt}(t)=\frac{\|\frac{d\psi}{dt}(t)\|}{c(1-\lambda)}\quad \mbox{if}\quad
t\in D_\lambda,\\
\frac{dS_\lambda}{dt}(t)=\frac{1}{(1-\lambda)}, \quad \mbox{if}\quad
t\in E_\lambda.
\ee
Clearly the map $S_\lambda:[0,1]\to [0,S_\lambda(1)]$ is one to one with
$S_\lambda(1)> 1$. Its inverse $T_\lambda:[0,S_\lambda(1)]\to [0,1]$.
For $s\in [0,S_\lambda(1)]$ we define
\be
\xi_\lambda(s)=\psi(T_\lambda(s)).
\ee
On the one hand
\be
\|\frac{d\xi_\lambda}{dt}(t)\|\leq \max(c(1-\lambda),\frac{1}{\lambda}).
\ee
On the other hand the hypothesis
\be
\int_0^1\tilde{L}(\psi(t),\frac{d\psi}{dt}(t))dt:=\sum_{i=1}^d
\int_0^1\tilde{L}_i(\psi(t),\frac{d\psi_i}{dt}(t))dt < +\infty
\ee
implies that for each $i$
\be
\int_0^1|\frac{d\psi_i}{dt}(t)-b_i^\eta(\psi(t))|dt < +\infty,
\ee
for $\eta=1,2$. From the triangular inequality
\be
\int_0^1|\frac{d\psi_i}{dt}(t)|dt\leq \int_0^1|\frac{d\psi_i}{dt}(t)-b_i^\eta(\psi(t))|dt+
\int_0^1|b_i^\eta(\psi(t))|dt
\ee
we derive for each $i$ that
\be
\int_0^1\|\frac{d\psi}{dt}(t)\|dt < +\infty.
\ee
Now the rest of the proof is the same as in Dupuis et al.
We prove
\be
\sup_{t\in [0,1]}\|\xi_\lambda(t)-\psi(t)\|\to 0
\ee
and
\be
\int_0^1\tilde{L}(\xi_\lambda(t),\frac{d\xi_\lambda}{dt}(t))dt\to
\int_0^1\tilde{L}(\psi(t),\frac{d\psi}{dt}(t))dt
\ee
as $\lambda\to 0$, which achieves the proof of Lemma 1.

{\bf Lemma 2.} Let $\xi\in C_{\x(0)}([0,1])$ such that $\int_0^1\tilde{L}(\xi(t),\frac{d\xi}{dt}(t))dt< +\infty$ and
\be
\sup_{t\in [0,1]}\|\frac{d\xi}{dt}(t)\| \leq \theta.
\ee
For any $\delta >0$ there exists $\sigma >0$ and $\varphi^\sigma\in \mathcal{N}_0$ such that
\be
\sup_{t\in [0,1]}|\xi(t)-\varphi^\sigma(t)|\leq \delta,\quad \mbox{and}\quad
\int_0^1\tilde{L}(\varphi^\sigma(t),\frac{d\varphi^\sigma}{dt}(t))dt\leq
\int_0^1\tilde{L}(\xi(t),\frac{d\xi}{dt}(t))dt+
2\delta
\ee

{\bf Proof.}

%We also recall that $\xi_i(t)=0$ implies $\xi_i'(t)=0$  $t$ a.s. Hence
%for $\xi_i(s)=0$ and for $|s-t|\leq \sigma_1$ implies
%\be
%|L_i^{0}(\xi(s),\xi'(s))-L_i^{0}(\xi(t),\xi'(s))|\leq \delta.
%\ee

We define for each $i$
\be
G_i^{0}&=&\{t\in [0,1]:\quad \xi_i(t)=0\},\\
G_i^1&=&\{t\in [0,1]:\quad \xi_i(t)<0\},\\
G_i^2&=&\{t\in [0,1]:\quad \xi_i(t)>0\}.
\ee
%\end{document}
Following Dupuis et al. for all $\sigma >0$ there exists $B_i^0=\bigcup_{j=1}^{J_i} [c_j(i),d_j(i)]$ such that
$\xi_i(c_j(i))=\xi_i(d_j(i))=0$
$d_j(i)-c_j(i)\leq \sigma$,
$d_j(i)\leq c_{j+1}(i)$. We suppose that $c_1(i)=0$ and $d_{J_i}(i)=1$.
We choose finitely many numbers $(e_j^k(i), k=1, \ldots, K_j(i))$
such that
\be
d_j(i)=e_j^1(i)< \ldots < e_j^{K_j(i)}(i)=c_{j+1}(i)
\ee
and
\be
e_j^{k+1}(i)-e_j^k(i) < \sigma.
\ee
%Note that for all $i$, we have $\xi_i(t)\neq 0$ for $t\in (d_j(i),c_{j+1}(i))$.

We define for each $i$
the function $\varphi_i^\sigma$ is piecewise linear interpolation of $\xi_i$ with interpolation points
\be
\{c_j(i),d_j(i),e_j^k(i):\quad j=1,\quad J_i, k=1, \ldots, K_j(i)\}.
\ee
For each $\delta >0$ there exists $\sigma_1>0$
such that for $0\leq \sigma\leq \sigma_1$ implies that
\be
\sup_{t\in [0,1]}|\xi_i(t)-\varphi_i^\sigma(t)|\leq \delta.
\ee
Clearly $\varphi_i^\sigma(t)=0$ for all $t\in [c_j(i),d_j(i)]$, and
$\varphi_i^\sigma(t)\neq 0$ for all $t\in [e_j^k(i),e_j^{k+1}(i)]$.

We define
\be
a_j^{\eta}(i)&=&\int_{c_j(i)}^{d_j(i)}{\bf 1}_{G_i^{\eta}}(t)dt,\\
\beta_j^\eta(i)&=&\frac{1}{a_j^\eta(i)}\int_{c_j(i)}^{d_j(i)}\frac{d\xi_i}{ds}(s){\bf 1}_{G_i^{\eta}}(s)ds.
\ee
Clearly
\be
\sum_{\eta=0,1,2}a_j^\eta(i)=d_j(i)-c_j(i).
\ee
Since $\xi_i(t)=0$ implies $\frac{d\xi_i}{dt}(t)=0$ a.s. then
\be
\beta_j^0(i)=0.
\ee
Since $G_i^1$, $G_i^2$ are open, then it can be written as a countable union of open
intervals at each endpoint of which
$\xi_i=0$. Hence $\beta_j^\eta(i)=0$ for $\eta=1,2$. It follows that
\be
L_i^1(\xi(c_j(i)),0)\geq L_i^0(\xi(c_j(i)),0),\\
L_i^2(\xi(c_j(i)),0)\geq L_i^0(\xi(c_j(i)),0).
\ee
We have for each $i,j$ that
\be
\int_{c_j(i)}^{d_j(i)}\tilde{L}_i(\xi(t),\frac{d\xi}{dt}(t))dt=\sum_{\eta=0,1,2}
\int_{c_j}^{d_j}{\bf 1}_{G_i^\eta}(t)L_i^\eta(\xi(t),\frac{d\xi_i}{dt}(t))dt.
\ee

From the continuity of $L_i^\eta$ and $\xi$ and
the fact that $\sup_{t\in [0,1]}\|\frac{d\xi}{dt}(t)\|\leq \theta$
there exists $\sigma_2\leq \sigma_1$
such that $\sigma\leq \sigma_2$ implies
\be
\sup_{t\in [0,1]} |L_i^\eta(\xi(t),\frac{d\xi_i}{dt}(t))- L_i^\eta(\xi(c_j(i)),\frac{d\xi_i}{dt}(t))|\leq \delta.
 \ee
 It follows that
\be
\int_{c_j(i)}^{d_j(i)}{\bf 1}_{G_i^\eta}(t)L_i^\eta(\xi(t),\frac{d\xi_i}{dt}(t))dt\geq
\int_{c_j(i)}^{d_j(i)}{\bf 1}_{G_i^\eta}(t)L_i^\eta(\xi(c_j(i),\frac{d\xi_i}{dt}(t))dt-(d_j(i)-c_j(i))\delta.
\ee

From the convexity of the function  $\beta_i\to L_i^\eta(\x,\beta_i)$
for each $\x$ fixed and for each $\eta=1,2$, we have
\be
\frac{1}{a_j^\eta(i)}\int_{c_j(i)}^{d_j(i)}{\bf 1}_{G_i^\eta}(t)L_i^\eta(\xi(c_j(i),\frac{d\xi_i}{dt}(t))dt
\geq L_i^\eta(\xi(c_j(i)),\beta_j^\eta(i))=L_i^\eta(\xi(c_j(i)),0)\geq L_i^0(\xi(c_j(i)),0).
\ee
For $\eta=0$, we have
\be
\int_{c_j(i)}^{d_j(i)}{\bf 1}_{G_i^0}(t)L_i^0(\xi(c_j(i),0)dt=
a_j^0(i)L_i^0(\xi(c_j(i),0).
\ee
Finally, we have
\ben
\int_{c_j(i)}^{d_j(i)}\tilde{L}_i(\xi(t),\frac{d\xi}{dt}(t))dt\geq
(d_j(i)-c_j(i))L_i^0(\xi(c_j(i)),0)-\delta(d_j(i)-c_j(i)).
\label{Lzeroxicji}
\een
Observe that
\be
\xi_i(c_j(i))=\varphi_i^\sigma(c_j(i)),
\ee
but there is no guaranty that
\be
\xi_l(c_j(i))=\varphi_l^\sigma(c_j(i))
\ee
for $l\neq i$. However there exist
$\alpha_l\leq c_j(i)\leq \beta_l$ such that
\be
\xi_l(\alpha_l)=\varphi_l^\sigma(\alpha_l),\quad
\xi_l(\beta_l)=\varphi_l^\sigma(\beta_l)\\
\beta_l-\alpha_l\leq \sigma.
\ee
More precisely
\be
\alpha_l=c_{j_l}(l),\quad \beta_l=d_{j_l}(l)
\ee
for some $j_l$ or
\be
\alpha_l=e_{j_l}^{k_l}(l),\quad \beta_l=e_{j_l}^{k_l+1}(l).
\ee
It follows that for small $\sigma$ we have
\be
\xi_l(c_j(i))\approx\varphi^\sigma(c_j(i))\approx
\varphi^\sigma(t),\quad \forall\,t\in [c_j(i),d_j(i)].
\ee
Then for small $\sigma$ we have
\be
 L_i^0(\xi(c_j(i)),0)\geq L_i^0(\varphi^\sigma(t),0)-\delta,\quad \forall\,t\in [c_j(i),d_j(i)].
\ee
Now the inequality (\ref{Lzeroxicji}) becomes
\be
\int_{c_j(i)}^{d_j(i)}\tilde{L}_i(\xi(t),\frac{d\xi}{dt}(t))dt\geq
\int_{c_j(i)}^{d_j(i)}L_i^0(\varphi^\sigma(t),\frac{d\varphi_i^\sigma}{dt}(t))dt-2\delta(d_j(i)-c_j(i))\\
=\int_{c_j(i)}^{d_j(i)}\tilde{L}_i(\varphi^\sigma(t),\frac{d\varphi_i^\sigma}{dt}(t))dt-2\delta(d_j(i)-c_j(i)).\\
\ee The same proof shows that \be
\int_{e_j^{k}(i)}^{e_j^{k+1}(i)}\tilde{L}_i(\xi(t),\frac{d\xi}{dt}(t))dt\geq
\int_{e_j^{k}(i)}^{e_j^{k+1}}
\tilde{L}_i(\varphi^\sigma(t),\frac{d\varphi_i^\sigma}{dt}(t))dt-2\delta(e_j^{k+1}(i)-e_j^k(i)).
\ee Finally \be &&\int_0^1\tilde{L}(\xi(t),\frac{d\xi}{dt}(t))dt=\sum_{i=1}^d
\int_0^1\tilde{L}_i(\xi(t),\frac{d\xi}{dt}(t))dt\\
&=&\sum_{i=1}^d\{\sum_{j=1}^{J_i}
(\int_{c_j(i)}^{d_j(i)}\tilde{L}_i(\xi(t),\frac{d\xi}{dt}(t))dt+
\sum_{k=1}^{K_j(i)}
\int_{e_j^k(i)}^{e_j^{k+1}(i)}\tilde{L}_i(\xi(t),\frac{d\xi}{dt}(t))dt)\}\\
&&\geq \int_0^1\tilde{L}(\varphi^\sigma(t),\frac{d\varphi^\sigma}{dt}(t))dt-2\delta.
\ee

Now back to the upper bound.
Let $\tau>0$ and
$\psi\in C_{\x(0)}([0,1])$ such that
\be
I_{\x(0)}(\psi) \leq \inf\{I_{\x(0)}(\varphi):\varphi\in C([0,1])\}+\tau,
\ee
and then
\be
\limsup_{\eps\to 0} H^\eps(\x(0))\leq \inf\{I_{\x(0)}(\varphi)+h(\varphi):\quad
\varphi\in C([0,1])\}+2\tau
\ee
for any $\tau>0$.

Finally we obtain for any continuous and bounded map $h$ that
\ben
 \lim_{\eps\to 0} H^\eps(\x(0))=\inf\{I_{\x(0)}(\varphi)+h(\varphi):\quad
\varphi\in C([0,1])\}.
\label{laplaceldp}
\een

\section{Different large deviation formulation}
In general, a family of probability measure $(\IP^\eps, \eps >0)$ on
a metric space $(X,\Delta)$ satisfies the large deviation principle
(LDP) with the rate function $I$ if the following conditions are
satisfied:

a) $I: X\to [0,+\infty)$ is lower semi-continuous,

b) For each $r>0$, $\{x\in X; I(x)\leq r\}$ is precompact,

c) For any $R>0$, there exists a compact set $K$ such that for any
$\delta >0$, we have  for small $\eps$,
\be
\IP^\eps(B_\delta^c(K))\leq \exp(-\frac{R^2}{\eps^2}),
\ee

d) $\lim_{\delta\to 0}\liminf_{\eps\to 0}\eps^2\ln\{\IP^\eps(B_\delta(x))\}=
\lim_{\delta\to 0}\limsup_{\eps\to 0}\eps^2\ln\{\IP^\eps(B_\delta(x))\}=-I(x)$.

Here $B_\delta(K)$ denotes the $\delta$-neighborhoods of any compact set $K$,
and $B_\delta^c(K)$ its complement.

Let $\IP_{\x(0)}^\eps$ be the probability distribution of $\x^\eps$
(\ref{lassodiffusion}) starting from $\x(0)$. It's known that the
limit (\ref{laplaceldp}) is equivalent to say that the family
$(\IP_{\x(0)}^\eps: \eps>0)$ satisfies the LDP with the rate
function $I_{\x(0)}$ (\ref{ratefunction}). See \cite{DZ} for a
general theory of the large deviation principle.

\end{document}